\documentclass[11pt,oneside]{article}

\newcommand{\documentdate}{8 February 2019}

\usepackage{a4wide,latexsym,amsmath,amssymb}

\newcommand{\numsection}[1]{\section{#1}\setcounter{equation}{0}}

\newtheorem{theorem}{Theorem}[section]
\newtheorem{definition}{Definition}
\newtheorem{proposition}[theorem]{Proposition}
\newtheorem{lemma}[theorem]{Lemma}
\newtheorem{corollary}[theorem]{Corollary}
\newtheorem{example}{Example}[section]
\newcommand{\beqn}[1]{\begin{equation}\label{#1}}
\newcommand{\eeqn}{\end{equation}}
\newcommand{\calA}{{\cal A}}
\newcommand{\calB}{{\cal B}}

\newcommand{\calD}{{\cal D}}
 
\newcommand{\calF}{{\cal F}} 
 
\newcommand{\calI}{{\cal I}} 
\newcommand{\calO}{{\cal O}} 
\newcommand{\calU}{{\cal U}}
\newcommand{\calV}{{\cal V}}

\newcommand{\calX}{{\cal X}}
\newcommand{\calY}{{\cal Y}}
\newcommand{\barf}{\overline{f}}
\newcommand{\req}[1]{(\ref{#1})}
\newcommand{\tim}[1]{\;\; \mbox{#1} \;\;}
\renewcommand{\Re}{\hbox{I\hskip -2pt R}}

\newcommand{\ii}[1]{\{1, \ldots, #1 \}}
\newcommand{\tr}{{\rm tr}}
\newcommand{\ms}{\;\;\;\;}
\def\E{\mathbb{E}}
\def\P{\mathbb{P}}
\def\S{\mathbb{S}}
\def\R{\mathbb{R}}
\def\T{\mathbb{T}}
\newcommand{\qed}{\hphantom{.}\hfill Q.E.D.\medbreak}
\newcommand{\ep}{{\,\Box\,}}
\newcommand{\epb}{\raisebox{-1pt}{$\Box$}}
\newcommand{\bep}{\stackrel{-}{\ep}}
\newcommand{\uep}{\,\underline{\Box}\,}

\everymath{\displaystyle}

\title{Bernstein Concentration Inequalities for Tensors\protect\\
       via Einstein Products}

\author{
Ziyan Luo\thanks{
  State Key Laboratory of Rail Traffic Control and Safety,
  Beijing Jiaotong University, Beijing, China.
  E-mail: zyluo@bjtu.edu.cn},
~Liqun Qi\thanks{
  Department of Applied Mathematics,
  The Hong Kong Polytechnic University,
  Hung Hom, Kowloon, Hong Kong.
  E-mail: liqun.qi@polyu.edu.hk},
~and~Philippe L. Toint\thanks{
  Namur Center for Complex Systems (naXys),
  University of Namur, 61, rue de Bruxelles, B-5000 Namur, Belgium.
  Email: philippe.toint@unamur.be}
}

\date{\documentdate}

\begin{document}
\maketitle
\begin{abstract}
A generalization of the Bernstein matrix concentration inequality to random tensors
of general order is proposed.  This generalization is based on the use of
Einstein products between tensors, from which a strong link can be established
between matrices and tensors, in turn allowing exploitation of existing
results for the former.
\end{abstract}

{\bf AMS subject classifications.} 15A52, 15A72, 49J55, 60H25.
{\bf Keywords}: random tensors, concentration inequality, Einstein products,
subsampling, computational statistics.

\section{Introduction}

The theory of random matrices has a rich history starting with Hurwitz (see
\cite{Forr10}) and Wishart \cite{Wish28} in the first half of the 20th
century. While it has developped on its own right within probability theory, it has
also found applications in many diverse domains of computational
statistics, ranging from matrix approximation \cite{FrieKannVemp98} to compressed sensing
\cite{Dono06}, graph theory \cite{AhlsWint02}, sparsification \cite{AchlMcSh01}
or subsampling of data \cite{WillSeeg01}.  Important tools in several of these fields
are matrix concentration theorems that give results on expectation, norm
distribution and probability of deviation from the expectation.  We refer the
interested reader to the excellent book by Tropp \cite{Trop15} for further
elaboration and an extensive bibliography.

The purpose of this short paper is to extend one of the proeminent matrix
concentration results, the Bernstein inequality, to the case of
tensors of general order.  This extension was originally motivated by the
desire to extend the use of the Bernstein inequality in subsampling estimation
of gradients and Hessians of additive multivariate real functions
\cite{BellGuriMori18,BellGuriMoriToin18,BellKrejKrkl18,ChenJianLinZhan18,
      KohlLucc17,XuRoosMaho17,XuRoosMaho18} to derivatives of
higher degree, thereby providing estimation tools for general Taylor's
expansions of such functions.  It is however clear that applications of the
new tensor result has wider potential, including, for instance, randomized
tensor sparsification (such as in video streaming) or randomized tensor
products for fast computations.

Our approaches hinges on Einstein products of tensors and associated
``matricization'' transformations: these recast tensors in the form of large
matrices to which known results of matrix concentration inequalities \cite{Trop15}
may then be applied.

The paper is organized as follows. Section~\ref{defs-s} introduces the
Einstein products and states some of its properties that are central to our
development. We then state the Bernstein concentration inequality for
Einstein-symmetric tensor of even order in Section~\ref{sBern-s}.  The more
general inequality for Einstein-symmetric tensors of arbitrary order is
derived in Section~\ref{gBern-s} and an ``intrinsic dimension'' version of
this inequality presented in Section~\ref{idBern-s}. Some conclusions and
perspectives are finally presented in Section~\ref{concl-s}.

\numsection{Tensors and the Einstein Product}\label{defs-s}

We start by defining the Einstein tensor product for high-order tensors, first
introduced by Lord Kelvin in 1856 \cite{Kelv56} and named after Albert
Einstein for his work in \cite{Eins07}.

\begin{definition}[Einstein Product, \cite{Eins07}]\label{Eproduct}
Let $\calA$ be a tensor in $\Re^{I_1\times\cdots\times I_m\times K_1\times\cdots\times K_m}$
and $\calB$ be a tensor in $\Re^{K_1\times\cdots\times K_m\times J_1\times\cdots\times J_p}$.
The Einstein product of $\calA$ and $\calB$, denoted by $\calA\ep\calB$,
is defined by
\beqn{EP}
\left(\calA\ep\calB\right)_{i_1\ldots i_mj_1\ldots j_p}
= \sum_{k_1\ldots k_m} a_{i_1\ldots i_m k_1\ldots k_m}b_{k_1\ldots k_mj_1\ldots j_p},
\tim{ for all } i_1,\ldots,i_m, j_1,\ldots,j_p,
\eeqn
\end{definition}
In this definition, each lowercase index varies from 1 to its uppercase
equivalent: for instance $i_2$ varies from $1$ to $I_2$, $k_3$ from $1$ to
$K_3$ and $j_1$ from $1$ to $J_1$.

The Einstein product can be regarded as a higher order generalization of the
standard matrix multiplication in which $m=p=1$. Such a contraction product
has been widely used in the areas of continuum mechanics \cite{LaiRubiKrem09}
and relativity theory \cite{Eins07}. Notice that in $\T_{m,d}$, the space of real
tensors of order $m$ and dimension $d$, that is the set of multiarrays $\calA
= (a_{i_1, \ldots, i_m})$ where $i_j$ varies from $1$ to $d$ for $j= 1,\ldots,
m$, the Einstein product satisfies the closure property
\[
\calA, \calB\in \T_{2m,d} \Longrightarrow \calA\ep\calB\in \T_{2m,d}.
\]
This nice property allows us to follow \cite{BrazLiNavaTamo13} and define
several new concepts based on the Einstein product for tensors.  
\begin{definition}\label{concepts} Let
$\calA=\left(a_{i_1\ldots i_mj_1\ldots j_m}\right)\in \T_{2m,n}$.
\begin{itemize}
  \item[(i)]{\bf Transpose:} The transpose of $\calA$, denotes by $\calA^\top$, is
    defined by the relations
    \[
    \left(\calA^\top\right)_{i_1\ldots i_mj_1\ldots j_m}
    =\left(\calA\right)_{j_1\ldots j_mi_1\ldots i_m}
    \tim{ for all } i_1,\ldots, i_m, j_1, \ldots, j_m.
    \]
  \item[(ii)]{\bf Einstein-Symmetric Tensor:} $\calA$ is called Einstein-symmetric, or
    \epb-symmetric, if $\calA^\top = \calA$.
    The set of all \epb-symmetric tensors in $\T_{2m,d}$ is a subspace and is denoted by $\S_{2m,d}$.
  \item[(iii)]{\bf Diagonal Tensor:} An \epb-symmetric tensor $\calA$ is
    said to be diagonal if $a_{i_1\ldots i_mj_1\ldots j_m}=0$ whenever
    $\prod_k\delta_{i_kj_k}=0$, where $\delta_{ij}$ is the Kronecker delta.
  \item[(iii)]{\bf Identity Tensor:} The Einstein-identity tensor, denoted by
    $\calI^\ep$, is a diagonal \epb-symmetric tensor with
    $a_{i_1\ldots i_mi_1\ldots i_m}=1$ for all $i_1, \ldots, i_m$. 
  \item[(iv)]{\bf Orthogonal Tensor:} $\calA$ is called Einstein-orthogonal,
    or \epb-orthogonal, if $\calA^\top\ep\calA = \calI^\ep$.
\item[(vi)]{\bf EVD:} If $\calA \in \S_{2m,d}$, then 
  \begin{equation}\label{EVD}
  \calA = \calU \ep \calD \ep \calU^\top
  \end{equation}
  is called an eigenvalue decomposition (EVD) of $\calA$, where $\calU$ is
  \epb-orthogonal and $\calD$ is \epb-symmetric and diagonal.
  Each $d_{i_1\ldots i_mi_1\ldots i_m}$ in $\calD$ is called an Einstein eigenvalue of
  $\calA$, or \epb-eigenvalue. The \epb-eigenvalues of $\calA$ are denoted by
  $\lambda^\ep_i(\calA)$ ($i\in\ii{d^m}$). 
\item[(vii)]{\bf Spectral norm and trace:} The Einstein-spectral norm and trace
  of $\calA \in \S_{2m,d}$ are defined by
  \[
  \|\calA\|^\ep = \max_{i\in \ii{d^m}} \left|\lambda^\ep_i(\calA)\right|
  \tim{ and }
  \tr^\ep(\calA) = \sum_{i\in \ii{d^m}}\lambda^\ep_i(\calA).
  \]
\end{itemize}
\end{definition}

As in \cite{BrazLiNavaTamo13}, we introduce the important bijective
``matricization'' transformation $f$ that maps each tensor $\calA\in
\T_{2m,d}$ to a matrix $A \in \R^{d^m\times d^m}$ with
$A_{ij} = a_{i_1\ldots i_m j_1\ldots j_m}$, where
\beqn{resp}
i = i_1 + \sum_{k=2}^m \left((i_k-1)d^{k-1}\right)
\tim{ and }
j = j_1 + \sum_{k=2}^m \left((j_k-1)d^{k-1}\right).
\eeqn
Note that
\beqn{f-of-vect}
f(x^{\otimes m})
= f( \underbrace{x \otimes \cdots \otimes x}_{m~{\rm times}})
= \underbrace{x \bullet \cdots \bullet x}_{m~{\rm times}}
= x^{\bullet m}
\tim{ for } x \in \Re^d,
\eeqn
where $\otimes$ denotes the tensor external product and $\bullet$ the
Kronecker product. Importantly for our purposes, it is proved in \cite{BrazLiNavaTamo13}
that 
\beqn{transform}
f(\calA\ep\calB) = f(\calA) \cdot f(\calB),
\eeqn
where $\cdot$ is the standard matrix multiplication. Thus the consistency of
the concepts introduced in Definitions~\ref{concepts} results from standard
matrix analysis.

The property \req{transform} in turn implies the following useful results.

\begin{proposition}\label{observations} Let
  $\calA=\left(a_{i_1\ldots i_mj_1\ldots j_m}\right)\in \T_{2m,d}$ be
  an \epb-symmetric tensor with EVD given by
  $\calA = \calU \ep\calD \ep \calU^\top$. We then have that
\begin{itemize}
  \item[(i)]  $f(\calA^\top) = f(\calA)^\top$, and hence
              $f(\calA) = f(\calU)\cdot f(\calD)\cdot f(\calU)^\top$;
  \item[(ii)] All eigenvalues of $f(\calA)$ are \epb-eigenvalues of
              $\calA$ and vice-versa;
  \item[(iii)] $\tr^\ep(\calA)= \sum_{i\in \ii{d^m}}\lambda_i(f(\calA))$;
  \item[(iv)] $\calA\ep\calA = \calU\ep \left(\calD\circ \calD\right)\ep\calU^\top$,
              where $\circ$ denotes the Hadamard product;
\end{itemize}
\end{proposition}

\noindent
Moreover, we may also establish a relation between
the Einstein- and the standard Z-eigenvalues. Here we simply recall that a real
scalar $\lambda$ is called a Z-eigenvalue of a symmetric real tensor
$\calA \in \T_{2m,d}$, if there exists real unit vector $x\in \R^d$ such that
\[
\calA x^{2m-1} = \lambda x,
\tim{ where } \calA x^{2m-1}
= \left(\sum_{i_2\ldots i_{2m}} a_{ii_2\ldots i_{2m}}x_{i_2}\cdots x_{i_{2m}}\right)\in\R^d
\]
(see \cite{QiLuo17}). As pointed out in \cite{Qi05}, Z-eigenvalues of
even-order symmetric real tensors always exist.

\begin{lemma}\label{bound} For an \epb-symmetric real tensor $\calA\in
  \T_{2m,d}$, we have that, whenever Z-eigenvalues of $\calA$ exist,
    $\lambda^\ep_{\max}(\calA)\geq \lambda^Z_{\max}(\calA)$.
\end{lemma}
\noindent{\emph{Proof.}}
By direct calculation, we have that
\begin{eqnarray}
  \lambda^\ep_{\max}(\calA)
    &   =  & \max_{y\in\R^{d^m}\setminus \{0\}}\frac{y^\top f(\calA)y}{\|y\|_2^2}
             \nonumber \\
    & \geq & \max_{x\in\R^d\setminus \{0\}} \frac{\langle f(\calA),
             \left(x^{\otimes m}\right)\cdot\left(x^{\otimes m}\right)^\top
             \rangle}{\|x^{\otimes m}\|_2^2 }
             \nonumber\\
    & \geq & \max_{x\in\R^d,~x^\top x =1} \calA x^{2m}\nonumber \\
    & \geq & \lambda_{\max}^Z(\calA), \nonumber
\end{eqnarray}
where the second inequality results from the observation that
$x^\top x = 1$ implies that $\|x^{\otimes m}\|_2^2 = 1$.
\qed

\numsection{The Bernstein Inequality for Even-Order Tensors}\label{sBern-s}

We now turn to random tensors, which are defined as follows. Let $(\Omega,
\calF, \P)$ be a probability space. A real $(m,d)$ random tensor $\calX$ is a
measurable map from $\Omega$ to $\T_{m,d}$.  A finite sequence $\{\calX_k\}$
of random tensors is independent whenever
\[
\P(\calX_k \in \calF_k \tim{ for all } k ) = \prod_k \P(\calX_k \in \calF_k)
\]
for every collection $\{\calF_k\}$ of Borel subsets of
$\T_{m,d}$. $\E(\calX)$, the expectation of the random tensor $\calX$, is, as
is the case for matrices, taken elementwise.

We are now in position to achieve our first objective: the Bernstein
inequality for even order real \epb-symmetric tensors based on Einstein products.  

\begin{theorem}\label{main} Consider a finite sequence $\{\calX_k\}$ of
independent random real \epb-symmetric tensors of order $2m$ and dimension $d$.
Assume that
\[
\E(\calX_k) = \calO
\tim{ and }
\lambda_{\max}^\ep(\calX_k)\leq L  \tim{for each } k.
\]
Consider the random tensor $\calY = \sum_k \calX_k$ and let $\nu(\calY)$ be
the tensor variance statistic of $\calY$ via Einstein product, that is
\[
\nu(\calY)
=\left\|\E(\calY^{\ep 2})\right\|^\ep
= \left\|\sum_k\E(\calX_k^{\ep 2})\right\|^\ep.
\]
Then
\beqn{expect}
\E\big(\lambda_{\max}^\ep(\calY)\big)
\leq\sqrt{2\nu(\calY)m\log d}+\frac{1}{3}Lm\log d.
\eeqn
Furthermore, for all $t\geq 0$,
\beqn{prob}
\P\left(\lambda_{\max}^\ep(\calY)\geq t\right)
\leq d^m\cdot \exp\left(\frac{-t^2/2}{\nu(\calY)+Lt/3}\right).
\eeqn
\end{theorem}
\noindent{\emph{Proof}.} First observe that the following equivalences
between tensors and matrices hold:
\begin{equation}\label{equvi}
  \E(\calX_k)=\calO \Longleftrightarrow \E(f(\calX_k)) = 0,
  \ms \|\calX_k\|^\ep\leq L \Longleftrightarrow \|f(\calX_k)\|\leq L,
  \ms \lambda_{\max}^\ep(\calY) = \lambda_{\max}(f(\calY)).
\end{equation}
Using those equivalences and applying the matrix Bernstein
inequality \cite[Theorem 6.6.1]{Trop15} to $f(\calY) = \sum_k f(\calX_k)$, we
then deduce the desired result. \qed

\noindent
Using Lemma~\ref{bound}, we then immediately deduce the following corollary
involving Z-eigenvalues.

\begin{corollary}\label{main-cor} Suppose that the assumptions of
  Theorem~\ref{main} hold and that Z-eigenvalues of $\calY$ exist. Then,
\beqn{expect-cor}
\E\big(\lambda_{\max}^Z(\calY)\big)
\leq\sqrt{2\nu(\calY)m\log d}+\frac{1}{3}Lm\log d.
\eeqn
Furthermore, for all $t\geq 0$,
\beqn{prob-cor}
\P\left(\lambda_{\max}^Z(\calY)\geq t\right)
\leq d^m\cdot \exp\left(\frac{-t^2/2}{\nu(\calY)+Lt/3}\right).
\eeqn
\end{corollary}

\noindent
This result reduces to the matrix Bernstein inequality for real symmetric
matrices \cite[Theorem 6.6.1]{Trop15} when $m=1$, since for any symmetric real
matrix $A$, $\calA$ is \epb-symmetric, 
\[
\lambda_{\max}^Z(A) =\lambda_{\max}^\ep(A) = \lambda_{\max}(A),
\ms
\|\E(A^{\ep 2})\|^\ep = \|\E (A^2)\|
\]
and $d^m=d$.

\numsection{The General Tensor Bernstein Inequality}\label{gBern-s}

As is the case for the matrix case, extending the condensation inequality to
tensors of odd order requires additional work. The notion of Einstein
product itself must first be extended to general tensors in $\T_{N,d}$.

\begin{definition}[Generalized Einstein Products]\label{generalized-Einstein}
  Let $\calA$, $\calB$ be two real tensors in $\T_{N,d}$, and
  $m = \left\lceil\frac{N}{2}\right\rceil$. Two generalized
  Einstein products of $\calA$ and $\calB$, denoted by $\calA \!\bep\! \calB$
  and $\calA \uep\calB$, are defined by
\beqn{gEinstein1}
\left(\calA \!\bep\! \calB\right)_{i_1\ldots i_m j_1\ldots j_m}
= \sum_{k_1,\ldots, k_{N-m}}
a_{i_1\ldots i_m k_1\ldots k_{N-m}}b_{j_1\ldots j_m k_1\ldots k_{N-m}}
\in \T_{2m,d},
\eeqn
and
\beqn{gEinstein2}
\left(\calA\uep \calB\right)_{k_1\ldots k_{N-m} k'_1\ldots k'_{N-m}}
= \sum_{i_1,\ldots, i_{m}}
a_{i_1\ldots i_m k_1\ldots k_{N-m}}b_{i_1\ldots i_m k'_1\ldots k'_{N-m}}
\in \T_{2(N-m),d},
\eeqn
respectively.
\end{definition}

\noindent
We examine two special cases.
\begin{itemize}
  \item[(i)] If $N=1$, the ranges between 1 and $N-m=0$
    in the above definition are interpreted as empty. In this case,
    $\calA$ and $\calB$ are vectors in $\Re^d$, say
    $a$ and $b$. Thus, $a \!\bep\! b = ab^\top$ and
    $a \uep b = a^\top b$, which are exactly the outer
    and inner products of vectors.
  \item[(ii)] If $N = 2m$, then $\calA$ and $\calB$ are in $\T_{2m,d}$, and
    \beqn{prodsok}
    \calA \!\bep\! \calB = \calA  \ep \calB^\top,
    \tim{ and }
    \calA \uep \calB = \calA^\top \ep \calB,
    \eeqn
    where $\calB^\top$ is defined in Definition~\ref{concepts} and $\ep$ is the
    Einstein product in Definition~\ref{Eproduct}, both for even-order tensors.
\end{itemize}

\noindent
We also need to generalize the bijective transformation $f$ which unfolds an
even-order tensor to a square matrix (as introduced in Section 2) to operate
on tensors of any order.  This is done as follows.

\begin{definition}[Matricization]\label{bijective2} Let $N \geq 1$, $d\geq 1$ and
  $m=\left\lceil \frac{N}{2}\right\rceil$. Define a bijective linear
  transformation $\barf$ from $T_{N,d}$  to $\Re^{d^m \times d^{N-m}}$
  such that for any tensor $\calA \in \T_{N,d}$,
  \[
  \left(\barf\left(\calA\right) \right)_{ik}
  = a_{i_1\ldots i_m k_1\ldots k_{N-m}},
  \]
  where
  \[
  i=i_1+\sum\limits_{l=2}^m \left((i_l-1)d^{l-1}\right)
  \tim{ and }
  k=k_{1}+\sum\limits_{l=2}^{N-m} \left((k_l-1)d^{l-1}\right).
  \]
\end{definition}

\noindent
Note that $\barf(\calA)$ need not be square or (obviously) symmetric.
As above, we consider two special cases.
\begin{itemize}
\item[(i)]  If $N=1$, the range between 1 and $N-m=0$
  is again interpreted as empty. It results that $\barf$
  is the identity transformation that maps any vector ${\bf x}\in\R^d$ to
  itself.
\item[(ii)] If $N=2m$, then $\barf$ coincides with the transformation $f$.
\end{itemize}

\noindent
The all important relation \req{transform} may also be generalized as follows.

\begin{lemma}\label{equiv} Let $N\geq 1$, $d\geq 1$ and $m=\left\lceil
  \frac{N}{2}\right\rceil$.  Then we have that, for all $\calA\in \T_{N,d}$,
\begin{equation}\label{eq}
f\left(\calA \!\bep\! \calA\right)
  = \barf\left(\calA\right)\cdot\barf\left(\calA\right)^\top
\tim{ and }
f\left(\calA \uep \calA\right)
  = \barf\left(\calA\right)^\top\cdot \barf\left(\calA\right).
\end{equation}
\end{lemma}
\noindent{\emph{Proof.}}  Notice that
\begin{eqnarray}
  \left(\calA\!\bep\! \calA\right)_{i_1\ldots i_m j_1\ldots j_m}
   &=& \sum_{k_1,\ldots, k_{N-m}} a_{i_1\ldots i_m k_1\ldots k_{N-m}}a_{j_1\ldots j_m k_1\ldots k_{N-m}}
       \nonumber \\
   &=& \sum_{k_1,\ldots, k_{N-m}} a_{j_1\ldots j_m k_1\ldots k_{N-m}}a_{i_1\ldots i_m k_1\ldots k_{N-m}}
        \nonumber\\
   &=& \left(\calA\!\bep\! \calA\right)_{j_1\ldots j_m i_1\ldots i_m}, \nonumber
\end{eqnarray}
for any $i_1, \ldots, i_m, j_1, \ldots, j_m$. Thus, $\calA\!\bep\! \calA\in \S_{2m,d}$
and hence $f\left(\calA\!\bep\! \calA\right)$ is well-defined.
Denote $B = f\left(\calA \!\bep\! \calA\right)$ and 
$C = \barf\left(\calA\right)\cdot \barf\left(\calA\right)^\top$.
From the definitions of $f$ and $\barf$, we know that the matrices
$B$ and $C$ have the same size, which is $d^m \times d^m$. For any $i$ and
$j \in \ii{d^m}$, there exist two $m$-tuples of indices
$\left(i_1,\ldots, i_m\right)$ and $\left(j_1,\ldots, j_m\right)$ that
uniquely determine by $i$ and $j$ via \req{resp}. By direct calculation, we
then obtain that
\begin{eqnarray}
C_{ij}
& = & \sum_{l=1}^{d^{N-m}}\left[\barf(\calA)\right]_{il}\left[\barf(\calA)\right]_{jl}
      \nonumber \\
& = & \sum_{k_1,\ldots, k_{N-m}} a_{i_1\ldots i_m k_1\ldots k_{N-m}}a_{j_1\ldots j_m k_1\ldots k_{N-m}}
      \nonumber\\
& = & \left(\calA\!\bep\! \calA\right)_{i_1\ldots i_m j_1\ldots j_m} \nonumber\\
& = & \left[f\left(\calA\!\bep\! \calA\right)\right]_{ij} \nonumber\\
& = & B_{ij}.
\end{eqnarray}
The proof for the case involving $\uep$ is similar. \qed

\noindent
We next need to revisit the definition of the spectral norm.

\begin{definition} Let $N, d\geq 1$ and
  $m=\left\lceil \frac{N}{2}\right\rceil$.
  Suppose that $\calA=\left(a_{i_1\ldots i_mk_1\ldots k_{N-m}}\right)\in \T_{N,d}$.
  The spectral norm of $\calA$ in the sense of generalized Einstein products,
  denoted as $\|\calA\|^{\bep}$, is defined by
  $\|\calA\|^{\bep} = \sqrt{\|\calA \!\bep\! \calA\|^\ep}$.
\end{definition}

\noindent
Using Proposition~\ref{observations} (iii) and \req{prodsok}, one verifies that
$\|\calA\|^{\bep}=\|\calA\|^{E}$ whenever $\calA \in \T_{2m,d}$.

As for the matrix case \cite{Trop15}, we now use a construct to build a
symmetric even-order object from (possibly) odd-order non-square parts.  This
is achieved by using the Hermitian dilation defined, for any real matrix $B$,
by
\[
H(B) = \left(\begin{array}{cc}
        O & B \\
        B^\top & O \\
\end{array}\right).
\]
It is then possible to establish a link between this construct and the
spectral norm just defined: first note that
\beqn{T2.1.28}
\lambda_{\max}(H(B)) = \|H(B)\| = \|B\|.
\eeqn
We may then use this identity to establish the following result.

\begin{lemma}
  Let $N, d\geq 1$ and $m=\left\lceil \frac{N}{2}\right\rceil$.
  Suppose $\calA=\left(a_{i_1\ldots i_mk_1\ldots k_{N-m}}\right)\in \T_{N,d}$.
  We then have that
  \begin{equation}\label{sym}
    \|\calA\|^{\bep}
    = \sqrt{\|\calA \uep \calA\|^\ep}
    = \|\barf\left(\calA\right)\|
    = \|H(\barf\left(\calA\right))\|
    = \lambda_{\max}\left(H\left(\barf\left(\calA\right)\right)\right).
  \end{equation}
\end{lemma}
\noindent{\emph{Proof.}} By direct calculation, we have
\begin{eqnarray}
\|\calA\|^{\bep} = \sqrt{\|\calA\!\bep\! \calA\|^\ep}
&=& \sqrt{\max_{i} \left|\lambda^\ep_i\left( \calA\!\bep\! \calA\right)\right|}
    \nonumber\\
&=& \sqrt{\max_{i} \left|\lambda_i\left(f\left( \calA\!\bep\! \calA\right)\right)\right|}
    \nonumber \\
&=& \sqrt{ \max_{i} \left|\lambda_i\left(\barf\left(\calA\right)\cdot
          \barf\left(\calA\right)^\top\right)\right| } \nonumber\\
&=&  \|\barf\left(\calA\right)\| \nonumber\\
&=&  \sqrt{\|\calA \uep \calA\|^\ep } \nonumber\\
\end{eqnarray}
where the second equality follows from Definition~\ref{concepts} (vii), the
third one by applying Proposition~\ref{observations} (ii), the fourth and the
sixth resulting from \eqref{eq}.  Now, using \req{T2.1.28},
\[
\|\barf\left(\calA\right)\|
= \|H(\barf\left(\calA\right))\|
= \lambda_{\max}\left(H\left(\barf\left(\calA\right)\right)\right),
\]
completing the proof. \qed

We are now in a position to state the general tensor
Bernstein inequality for random tensors of any order.

\begin{theorem}\label{main2}
Consider a finite sequence $\left\{\calX_k\right\}$ of independent
random tensors in $\T_{N,d}$ and let $m = \left\lceil \frac{N}{2}\right\rceil$.
Assume that, for some constant $L \geq 0$,
\[
\E(\calX_k) = \calO
\tim{ and }
\|\calX_k\|^{\bep} \leq L
\tim{ for all } k.
\]
Consider now the random tensor $\calY = \sum_k \calX_k$ and let $v(\calY)$ be the
generalized tensor variance statistic of the sum given by
\begin{eqnarray}
\nu(\calY)
& = & \max\left\{\left\|\E\left(\calY\!\bep\! \calY\right)\right\|^\ep,
                 \left\|\E\left(\calY\uep \calY\right)\right\|^\ep\right\}\nonumber\\
& = & \max\left\{\Big\|\sum_k\E\Big(\calX_k\!\bep\! \calX_k\Big)\Big\|^\ep,
                 \Big\|\sum_k\E\Big(\calX_k\uep \calX_k\Big)\Big\|^\ep\right\}\nonumber
\end{eqnarray}
Then
\begin{equation}\label{expect2}
  \E(\|\calY\|^{\bep})
  \leq\sqrt{2\nu(\calY)\log\left(d^m+d^{N-m}\right)}+\frac{1}{3}L\log\left(d^m+d^{N-m}\right).
\end{equation}
Furthermore, for all $t\geq 0$,
\begin{equation}\label{prob2}
  \P\left(\|\calY\|^{\bep}\geq t\right)
  \leq  \left(d^m+d^{N-m}\right)\cdot \exp\left(\frac{-t^2/2}{\nu(\calY)+Lt/3}\right).
\end{equation}
\end{theorem}
\noindent{\emph{Proof.}} The desired result follows from applying
\cite[Theorem~6.1.1]{Trop15} to the random matrix
$\barf(\calY) = \sum_k \barf(\calX_k)$ and using the facts that
\[
\|\calX_k\|^{\bep} = \|\barf(\calX_k)\|,
\ms
\|\calY\|^{\bep} = \|\barf(\calY)\|,
\]
and that
\begin{eqnarray}
\nu(\calY)
& = & \max\left\{\left\|\E\left(\barf(\calY)\cdot\barf(\calY)^\top\right)\right\|,
                \left\|\E\left(\barf(\calY)^\top\cdot \barf(\calY)\right)\right\|\right\}
                \nonumber\\
& = & \max\left\{\Big\|\E\Big(\sum_k\barf(\calX_k)\cdot\barf(\calX_k)^\top\Big)\Big\|,
                \Big\|\E\Big(\sum_k\barf(\calX_k)^\top\cdot \barf(\calX_k)\Big)\Big\|\right\}.
\end{eqnarray}
\qed

\noindent
Observe that the dimension-dependent factor on the right-hand side of
\req{prob} is $d^m+d^{N-m}$, which is larger than $md$, the factor one might
naively expect as a generalization of the matrix case, where this factor is
$2d$.  This larger bound somewhat limits the applicability of the results to
moderate values of $d$ and $m$. It is however worthwhile to note that we
have merely assumed the \epb-symmetry of the random tensors under
consideration, which is weaker than true symmetry.

\numsection{The Tensor Bernstein Inequality in Intrinsic Dimension}\label{idBern-s}

The above discussion about the dimension-dependent factor of \req{prob}
prompts the question of the extension of a version of the Bernstein
inequality where this factor can be improved.  This is the case of ``intrinsic
dimension'' version of this result, which we now consider.

Our approach first introduces Einstein-positive-(semi)definite tensors.
The positive semi-definiteness of real tensors has been discussed in
\cite{Qi05} and shown to have applications such as in biomedical imaging
\cite{QiYuWu10}. Recall that a real tensor $\calA\in \T_{2m,d}$ is called positive
semi-definite (PSD) if 
\[
\calA x^{2m} =
\sum_{i_1\ldots i_mj_1\ldots j_m}a_{i_1\ldots i_m j_1\ldots j_m}
x_{i_1}\cdots x_{i_m}x_{j_1}\cdots x_{j_m}
\geq 0, \tim{ for all } x\in\R^d.
\]
(see \cite{QiLuo17}).
Moreover, it has been shown in \cite{Qi05} that an even-order symmetric real
tensors is PSD if and only if all Z-eigenvalues are nonnegative.
Similarly, we can define such a nonnegativity in the sense of Einstein
products as follows.
\begin{definition}
An Einstein-symmetric tensor $\calA\in\T_{2m,d}$ is called Einstein-positive
semi-definite (\epb-PSD) (\epb-positive-definite (\epb-PD), respectively) if
and only if all its Einstein-eigenvalues are nonnegative (positive, respectively).
\end{definition}
We adopt the notation $\calA\succeq^\ep (\succ^\ep) \calO$ to represent that
$\calA$ is \epb-PSD (\epb-PD), and similarly $\calA\succeq^\ep (\succ^\ep) \calB$
is $\calA-\calB$ if \epb-PSD (\epb-PD). Such an \epb-PSD (\epb-PD) property is
actually stronger than the original PSD (PD) property, as stated in the
following lemma.

\begin{lemma}\label{psd}
Suppose that $\calA\in\S_{2m,d}$. If $\calA$ is \epb-PSD (\epb-PD), then $\calA$ is PSD (PD).
\end{lemma}
\noindent{\emph{Proof.}} Because of \req{transform}, $\calA$ is \epb-PSD if and only
if $f(\calA)$ is a PSD matrix. Then, for any $x\in\R^d$, it follows that 
\[
\calA x^{2m}
=  \langle f(\calA), \left(x^{\otimes m}\right)\left(x^{\otimes m}\right)^\top\rangle
= \left(x^{\otimes m}\right)^\top f(\calA) \left(x^{\otimes m}\right)
\geq  0.
\]
The proof for the \epb-PD case is similar. \qed

\noindent
Note that the PSD property does not, in general, imply the \epb-PSD property.
The following counterexample is taken from \cite[Example 4.5]{LuoQiYe15}.
\begin{example}\label{countereg}
Let $\calA\in \S_{4,3}$ with
$a_{1122}=a_{1212}=a_{1221}=a_{2112}=a_{2121}=a_{2211}=1$ and other entries
$0$. It is easy to verify that
$\calA x^{4} = 6x_1^2x_2^2\geq 0$ for any $x \in \R^3$,
whereas $y^\top f(\calA)y = 2y_1y_5<0$ for $y= (1,0,0,0,-1,0,0,0,0)^\top$.
\end{example}

\begin{proposition}\label{AA-l}
Let $\calA \in \T_{N,d}$ and $m=\left\lceil\frac{N}{2}\right\rceil$.  Then
$\calA \!\bep\! \calA \in \T_{2m,d}$ and 
$\calA \uep \calA \in \T_{2(N-d),d}$
are both \epb-PSD and PSD.
\end{proposition}
\noindent{\emph{Proof.}}
\[
\begin{array}{lcl}
(\calA \!\bep\! \calA) x^{2m}
&   =  & \langle f(\calA \!\bep\! \calA),
         \left(x^{\otimes m}\right)\left(x^{\otimes m}\right)^\top\rangle\\*[1.5ex]
&   =  & \left(x^{\otimes m}\right)^\top \barf(\calA)\cdot
         \barf(\calA)^\top \left(x^{\otimes m}\right)\\*[1.5ex]
&   =  & \| \barf(\calA)^\top \left(x^{\otimes m}\right)\|^2\\*[1.5ex]
& \geq & 0.
\end{array}
\]
The proof is similar for $\calA \uep \calA$.
\qed

Armed with these extended notions and the fundamental relation \req{transform}
applied to the Einstein EVD, we finally state an intrinsic-dimension version
of the Bernstein concentration inequality for tensors.

\begin{theorem}\label{main3}
Consider a finite sequence $\left\{\calX_k\right\}$ of independent
random tensors in $\T_{N,d}$ and let $m = \left\lceil \frac{N}{2}\right\rceil$.
Assume that, for some constant $L \geq 0$,
\[
\E(\calX_k) = \calO
\tim{ and }
\|\calX_k\|^{\bep} \leq L
\tim{ for all } k.
\]
Consider now the random tensor $\calY = \sum_k \calX_k$ and let $\calV_1$ and
$\calV_2$ be upper bounds for the tensor-valued variance statistics of $\calY$
introduced in Theorem~\ref{main2}, that is
\beqn{Vbounds}
\calV_1
\succeq^\ep \E\left(\calY\!\bep\! \calY\right)
= \sum_k\E\left(\calX_k\!\bep\! \calX_k\right) \tim{ and }
\calV_2
\succeq^\ep \E\left(\calY\uep \calY\right)
= \sum_k\E\left(\calX_k\uep \calX_k\right).
\eeqn
Let
\[
\nu(\calY) = \max\left\{\|\calV_1\|^\ep,\|\calV_2\|^\ep\right\}
\tim{ and }
d_{\calV}(\calY)
= \frac{1}{\nu(\calY)}
\Big(\tr^\ep(\calV_1)+\tr^\ep(\calV_2)\Big).
\]
Then, for $t \geq \sqrt{\nu(\calY)}+L/3$,
\begin{equation}\label{prob3}
  \P\left(\|\calY\|^{\bep}\geq t\right)
  \leq  4d_{\calV}(\calY)\cdot \exp\left(\frac{-t^2/2}{\nu(\calY)+Lt/3}\right).
\end{equation}
\end{theorem}
\noindent{\emph{Proof.}} We first observe that, because of
Proposition~\ref{AA-l}, $\E\left(\calY \!\bep\! \calY\right)$
and $\E\left(\calY \uep \calY\right)$ are
\epb-positive-semidefinite, which make the \epb-PSD ordering in \req{Vbounds}
well-defined. We also note that $d_{\calV}(\calY)$ is identical to the
intrinsic dimension of the matrix 
\beqn{V-def}
V = \left(\begin{array}{cc}
          \barf(\calV_1)^T &    0      \\
                0        & \barf(\calV_2)
          \end{array}
     \right),
\eeqn
where the standard (matrix) intrinsic dimension of a
positive-semidefinite matrix $M$ is the ratio $\tr(M)/\|M\|$.
The desired result then again follows from applying
an existing result for matrices (here \cite[Theorem~7.3.1]{Trop15})
to the random matrix $\barf(\calY) = \sum_k \barf(\calX_k)$.
\qed

\noindent
The main differerence between this theorem and Theorem~\ref{main2} is the
replacement of \req{prob2} by \req{prob3}: have to relax the range of $t$ for
which the inequality is valid but often gain in the ``dimension-dependent''
factor, since $d_{\calV}(\calY)$ never exceeds $d^m+d^{N-m}$ and can be much
smaller if $V$ in \req{V-def} is close to being of low rank.

\numsection{Conclusion}\label{concl-s}

We have considered the Einstein tensor products and reviewed the strong link
this concept establishes between standard matrix theory and tensor analysis.
This link has allowed us to restate the powerful Bernstein matrix
concentration inequality in the case of general tensors of arbitrary order.

Other concentration inequalities do exist for matrices (see \cite{Trop15} for
an overview). Whether they can be extended to tensors using a similar
approach, although likely, remains open at this stage.

It is interesting (and challenging) to examine if a better ``dimension
factor'' (closer to $md$) could be achieved by an approach where one does not
merely unfold tensors to matrices and use existing concentration results for these, but
where a true analysis of the tensor case is conducted.  The main difficulty
is to find an eigenvalue decomposition of (random) tensors with a a number of
``eigenvalues'' smaller than $d^m$ (this is for instance not necessarily the
case of Z-eigenvalues \cite{CartStur13}).

If one is to judge by the vast diversity of applications where matrix
concentration inequalities have been useful, our new result 
potentially opens several research paths in high-dimensional
computational statistics and numerical optimization.  In particular, its
application to sub-sampling methods for the estimation of
derivative tensors beyong the Hessian may now be considered, as it makes algorithms
based on high-order Taylor's expansions and models practical. The complexity
of optimization methods of this type has been analyzed in
\cite{BellGuriMoriToin18}, but the necessary probabilistic estimation properties were
so far limited to quadratic models.  The new tensor concentration inequality
thus allows further developements in a framework which is central to computational
deep learning.

{\footnotesize
  
\section*{Acknowledgment}

This research was partially supported by National Natural Science Foundation
of China (11771038, 11431002), the State Key Laboratory of Rail Traffic
Control and Safety, Beijing Jiaotong University (RCS2017ZJ001), and the Hong
Kong Research Grant Council (Grant No. PolyU 15300715, 15301716 and
15300717). The third author gratefully acknowledges the support of the
Hong Kong Polytechnic for the visit during which this research was initiated.


}

\end{document}